\documentclass[12pt]{amsart}

\usepackage{amsmath,amsthm, amssymb, amsfonts, epsfig,amscd}
\usepackage[all]{xy}

\usepackage[mathscr]{eucal}
\usepackage{mathrsfs}
\usepackage{latexsym}
\usepackage{graphicx}
\newdimen\harrowlength \harrowlength=60pt
\newdimen\varrowlength \varrowlength=.618\harrowlength
\newdimen\sarrowlength \sarrowlength=\harrowlength

\newdimen\hgrid \hgrid=15pt
\newdimen\vgrid \vgrid=15pt

\newdimen\hchannel  \hchannel=0pt
\newdimen\vchannel  \vchannel=0pt
\newdimen\channelwidth \channelwidth=3pt

\theoremstyle{plain}
\newtheorem{thm}[subsection]{Theorem}
\newtheorem{lem}[subsection]{Lemma}
\newtheorem{prop}[subsection]{Proposition}
\newtheorem{cor}[subsection]{Corollary}
\newtheorem{lemma}[subsection]{Lemma}
\newtheorem{theorem}[subsection]{Theorem}

\newtheorem*{conj10}{Conjecture of \cite{g-l}}

\theoremstyle{definition}

\newtheorem{defn}[subsection]{Definition}
\newtheorem{definition}[subsection]{Definition}
\newtheorem{remark}[subsection]{Remark}
\newtheorem{fact}[subsection]{Fact}

\def\L{\mathcal{L}}
\newcommand{\Z}{{\mathbb Z}}
\newcommand{\R}{{\mathbb R}}

\newcommand{\N}{{\mathbb N}}

\def\a{\alpha}
\def\t{\tau}

\def\f{\phi}
\def\l{\lambda}

\def\m{\mu}

\def\v{\vee}
\def\w{\wedge}

\def\ni{\noindent}

\begin{document}
\title[Singular loci of B-H toric varieties]{Singular loci of Bruhat-Hibi toric varieties}
\address{Department of Mathematics\\ Northeastern University\\ Boston, MA 02115}
\author[J. Brown]{J. Brown}
\author[V. Lakshmibai]{V. Lakshmibai${}^{\dag}$}

\email{brown.justin1@neu.edu, lakshmibai@neu.edu}
\thanks{${}^{\dag}$ Partially supported
by NSF grant DMS-0652386 and Northeastern University RSDF 07-08.}

\begin{abstract}
For the toric variety $X$ associated to the Bruhat poset of
Schubert varieties in a minuscule $G/P$, we describe the singular
locus in terms of the faces of the associated polyhedral cone.  We
further show that the singular locus is pure of codimension 3 in
$X$, and the generic singularities are of cone type.
\end{abstract}
\maketitle

\section*{Introduction}

Let $K$ denote the base field which we assume to be algebraically
closed of arbitrary characteristic. Given a distributive lattice
$\mathcal{L}$, let $X(\mathcal{L})$ denote the affine variety in
$\mathbb{A}^{\#\mathcal{L}}$ whose vanishing ideal is generated by
the binomials $X_\tau X_\varphi
-X_{\tau\vee\varphi}X_{\tau\w\varphi}$ in the polynomial algebra
$K[X_\alpha,\alpha\in\mathcal{L}]$ (here, ${\tau\vee\varphi}$
(resp. ${\tau\w\varphi}$) denotes the \emph{join} - the smallest
element of $\mathcal{L}$ greater than both $\tau,\varphi$ (resp.
the \emph{meet} - the largest element of $\mathcal{L}$ smaller
than both $\tau,\varphi$)). These varieties were extensively
studied by Hibi in  \cite{Hi} where Hibi proves that
$X(\mathcal{L})$ is a normal variety. On the other hand,
Eisenbud-Sturmfels show in \cite{ES} that a binomial prime ideal
is toric (here, ``toric ideal" is in the sense of \cite{St}). Thus
one obtains that $X(\mathcal{L})$ is a normal toric variety.  We
shall refer to such a $X\left(\mathcal{L}\right)$ as a \emph{Hibi
toric variety}.

For $\mathcal{L}$ being the Bruhat poset of Schubert varieties in
a minuscule $G/P$, it is shown in \cite{GLdef} that
$X(\mathcal{L})$ flatly deforms to ${\widehat{G/P}}$ (the cone
over $G/P$), i.e., there exists a flat family over
${\mathbb{A}}^1$ with ${\widehat{G/P}}$ as the generic fiber and
$X(\mathcal{L})$ as the special fiber. More generally, for a
Schubert variety $X(w)$ in a minuscule $G/P$, it is shown in
\cite{GLdef} that $X({\mathcal{L}}_w)$ flatly deforms to
${\widehat{X(w)}}$, the cone over $X(w)$ (here, ${\mathcal{L}}_w$
is the Bruhat poset of Schubert subvarieties of $X(w)$). In a
subsequent paper (cf. \cite{g-l}), the authors of loc.cit.,
studied the singularities of $X(\mathcal{L}), \mathcal{L}$ being
the Bruhat poset of Schubert varieties in the Grassmannian;
further, in loc.cit., the authors gave the following conjecture on
the singular locus of $X\left(\L\right)$:

\begin{conj10}\[\mbox{Sing}\, X\left(\L\right)=
\bigcup_{\left(\alpha,\beta\right)} Z_{\alpha,\beta} ,\] where
$\left(\alpha,\beta\right)$ is an (unordered) incomparable  pair
of join-meet irreducibles in $\L$, and $Z_{\alpha,\beta} = \{P\in
X\left(\L\right)\mid P\left( \theta\right) = 0,\, \forall\,
\theta\in [\alpha\wedge\beta, \alpha\vee\beta ]\}$.
\end{conj10}
(Here, for a $P\in X({\L})\subset \mathbb{A}^{\#{\L}}$, and
$\theta\in{\L}, P\left( \theta\right)$ denotes the $\theta$-th
co-ordinate of $P$.)

The sufficiency part of the above conjecture for the Bruhat poset
of Schubert varieties in the Grassmannian is proved in \cite{g-l},
using the Jacobian criterion for smoothness, while the
 necessary part of the conjecture is proved in \cite{Font}, using certain
desingularization of $X(\mathcal{L})$.

In \cite{G-H-paper}, the authors gave a simple proof of the above
conjecture for the Bruhat poset of Schubert varieties in the
Grassmannian using just the combinatorics of the polyhedral cone
associated to $X\left(\L\right)$.

It turns out that the above conjecture does not extend to a
general Hibi toric variety $X(\mathcal{L})$ (see \S 10 of \cite{G-H-paper} 
for a counter example). In \cite{G-H-paper}, the authors
conjectured that the above conjecture holds for other minuscule
posets. The main result of this paper is the proof of the above
conjecture for $\L$ being the Bruhat poset of Schubert varieties
in a minuscule $G/P$ (cf. Theorem \ref{theth}); we refer to the
corresponding $X(\mathcal{L})$ as a \emph{Bruhat-Hibi toric
variety}.  In fact, we show (cf. Theorem \ref{mainth}) that the
above conjecture holds for more general $X\left(\L\right)$,
namely, $\L$ being a distributive lattice such that
$J\left(\L\right)$ (the poset of join irreducibles) is a grid
lattice (see \S \ref{grid} for the definition of a grid lattice).
We further prove (cf. Theorem \ref{mainth}) that the singular
locus of $X(\mathcal{L})$ is pure of codimension 3 in
$X(\mathcal{L})$, and that the generic singularities are of cone
type (more precisely, the singularity type is same as that at the
vertex of the cone over the quadric surface $x_1x_4-x_2x_3=0$ in
$\mathbb{P}^3$ ).

\vskip.2cm\ni\textbf{Sketch of proof of the above conjecture for
Bruhat-Hibi toric varieties:} Let $\L$ be the distributive lattice
of Schubert varieties in a minuscule $G/P$, or more generally, a
distributive lattice such that the poset of join irreducibles is a
grid lattice. Let $T$ denote the torus acting on the toric variety
$X\left(\L\right)$. Let $M$ be the character group of $T$. Let
$\sigma $ be the polyhedral cone associated to the toric variety
$X\left(\L\right)$. If $\sigma^{\vee}$ is the cone dual to
$\sigma$ and $S_\sigma =\sigma^{\vee}\cap M$, then
$K[X\left(\L\right)]$ is the semigroup algebra $K[S_\sigma ]$. For
a face $\tau$ of $\sigma$, let $D_\tau =\{\alpha\in\L\mid
P_\tau\left(\alpha\right)\neq 0\}$, where $P_\tau$ (cf. $\S$
\ref{distinguished}) is the center of the orbit $O_\tau$. Now
$X_\tau$, the toric variety associated to the cone $\tau$, is open
in $X_\sigma$ ($=X\left(\L\right)$ ). Thus $X_\sigma$ is smooth at
$P_\tau$ if and only if $X_\tau$ is smooth at $P_\tau$; further,
$X_\tau$ is smooth at $P_\tau$ if and only if $X_\tau$ is
non-singular.

For $\tau$ such that $D_\tau = \L _{\alpha,\beta} = \L\setminus
[\alpha\wedge\beta , \alpha\vee\beta ]$, where $(\alpha,\beta)$ is an incomparable pair
of join-meet irreducibles in ${\L}$, we first determine a set of
generators for $\tau$ as a cone, and show that $X_\tau$ is a
singular variety. Conversely, if $\tau$ is such that $D_\tau$ is
not contained in any $\L_{\alpha,\beta}$, we show that $X_\tau$ is
non-singular. Thus the above conjecture is proved. As a
consequence, we obtain that Sing$\, X\left(\L\right)$ is pure of
codimension $3$ in $X\left(\L\right)$.

It should be remarked that the Hibi toric varieties are studied in
\cite{Wa} also where the author proves that the singular locus of
a Hibi toric variety has codimension at least three.

The sections are organized as follows: In \S \ref{first-sect}, we
recall some generalities on affine toric varieties. In \S
\ref{varlat}, we introduce the Hibi toric varieties, and recollect
some of the results (cf. \cite{LM}) on Hibi toric varieties
required for our discussion. In \S \ref{grid}, we introduce grid
lattices and prove some preliminary results on a distributive
lattice whose poset of join irreducibles is a grid lattice. In \S
\ref{faces}, we determine the singular locus of $X({\L}),{\L}$
being as above. In \S \ref{mini}, we apply the results of \S
\ref{faces} to Bruhat-Hibi toric varieties and determine the
singular loci of these varieties.

\vskip.2cm\ni \textbf{Acknowledgement}: The authors thank the
referee for many useful comments on the first version of this
paper.

\section{Generalities on toric varieties}\label{first-sect}

Since our main object of study is a certain affine toric variety,
we recall in this section some basic definitions on affine toric
varieties. Let $T=(K^*)^m$ be an $m$-dimensional torus.
\begin{defn}\label{equi} (cf. \cite{F}, \cite{KK})
An {\em equivariant affine embedding}  of a torus $T$ is an affine
variety $X\subseteq{\mathbb{A}}^l$ containing $T$ as a dense open
subset and equipped with a $T$-action $T\times X\to X$ extending
the action $T\times T\to T$ given by multiplication. If in
addition $X$ is normal, then $X$ is called an {\em affine toric
variety}.
\end{defn}

\subsection{The Cone Associated to a Toric Variety}\label{comb}
Let $M$ be the character group of $T$, and $N$ the $\Z$-dual of
$M$.  Recall (cf. \cite{F}, \cite{KK}) that there exists a
strongly convex rational polyhedral cone $\sigma \subset N_{\R}
(=N\otimes _{\Z}\R )$ such that
\[ K[X]=K[S_{\sigma}], \]
where $S_{\sigma}$ is the subsemigroup $\sigma^{\v}\cap M$,
$\sigma^{\v}$ being the cone in $M_{\R}$ dual to $\sigma$, namely,
$\sigma^{\v}=\{f\in M_{\R}\,|\,f(v)\ge 0, v\in
 \sigma\}$. Note that $S_{\sigma}$ is a finitely generated subsemigroup in
$M$.

\subsection{Orbit Decomposition in Affine Toric Varieties}\label{orbit}
We shall denote $X$ also by $X_\sigma$.  We may suppose, without
loss of generality, that $\sigma$ spans $N_{\R}$ so that the
dimension of $\sigma$ equals $\dim N_{\R}=\dim T$. (Here, by
dimension of $\sigma$, one means the vector space dimension of the
span of $\sigma$.)

\subsection{The distinguished point $P_\tau$}\label{distinguished}
  Each face $\tau$ determines a (closed) point
$P_\tau$ in $X_\sigma$, namely, it is the point corresponding to
the maximal ideal in $K[X](=K[S_\sigma])$ given by the kernel of
$e_\tau:K[S_\sigma]\rightarrow K $, where for $u\in S_\sigma$, we
have
$$e_\tau(u)=\begin{cases}1,&{\mathrm{\ if\ }}u\in \tau^{\perp}\\
0,&{\mathrm{otherwise}} \end{cases}$$ (here, $\tau^{\perp}$
denotes $\{u\in M_{\R}\,|\,u(v)=0,\forall v\in\tau\}$)

\subsection{Orbit Decomposition}\label{orbit-decomp} Let $O_\tau$ denote
the $T$-orbit in $X_\sigma$ through $P_\tau$. We have the
following orbit decomposition in $X_\sigma$:
$$\begin{gathered}
X_\sigma={\underset{\theta\le\sigma }{\cup}}\, O_\theta\\
{\overline{O_\tau}}={\underset{\theta\ge\tau }{\cup}}\, O_\theta\\
\dim\,\tau + \dim\,O_\tau=\dim\,X_\sigma
\end{gathered}$$
See \cite{F}, \cite{KK} for details.

\section{The toric variety associated to a distributive lattice}\label{varlat}
 We shall now study a special class of
toric varieties, namely, the toric varieties associated to
distributive lattices. We shall first collect some definitions as
well as some notation. Let $(\mathcal{L},\le)$ be a poset, i.e, a
finite partially ordered set. We shall suppose that $\mathcal{L}$
is {\em bounded}, i.e., it has a unique maximal, and a unique
minimal element, denoted $\widehat{1}$ and $\widehat{0}$
respectively. For $\mu,\lambda\in \mathcal{L},\mu\le\lambda$, we
shall denote
$$[\mu,\lambda]:=\{\tau\in\mathcal{L},\mu\le\tau\le\lambda\}$$ We
shall refer to $[\m,\l]$ as the\/ {\em interval from $\m$ to
$\l$}.
\begin{defn}\label{cover-def}
 The ordered  pair $(\l,\m)$ is called a\/ {\em cover}
 (and we also say
that $\l$\/ {\em  covers} $\m$ or $\mu$ is \emph{covered} by
$\lambda$) if $[\mu,\lambda]=\{\mu,\lambda\}$.
\end{defn}

\subsection{Distributive lattices}\label{s2}
\begin{defn} A {\em lattice\/} is a partially ordered set
$(\mathcal{L},\le)$ such that,
for every pair of elements $x,y\in \mathcal{L}$, there exist
elements $x\v y$ and $x\w y$, called the\/ {\em join},
respectively the\/ {\em meet} of $x$ and $y$, defined by:
\begin{gather}
x\v y\ge x,\ x\v y\ge y,\text{ and if  }z\ge x \text{ and }
z\ge y,\text{ then } z\ge x\v y,\notag\\
x\w y\le x,\ x\w y\le y,\text{ and if  }z\le x \text{ and } z\le
y,\text{ then } z\le x\w y.\notag
\end{gather}
\end{defn}

\begin{defn}\label{5.21}
Given a lattice ${\mathcal{L}}$, a subset ${\mathcal{L}}'\subset
{\mathcal{L}}$ is called a\/ {\em sublattice} of $\mathcal{L}$ if
$x,y\in{\mathcal{L}}'$ implies $x\w y\in{\mathcal{L}}'$, $x\v
y\in{\mathcal{L}}'$; $\L '$ is called an {\em embedded sublattice
of} $\L$ if
$$\t,\,\f\in\L,\quad\t\v\f,\,\t\w\f\in\L '\quad\Rightarrow\quad\t,
\,\f\in\L '.$$
\end{defn}

It is easy to check that the operations $\v$ and $\w$ are
commutative and associative.

\begin{defn}\label{dist-def} A lattice is called\/ {\em distributive} if the following
identities hold:
\begin{align}
x\w (y\v z)&=(x\w y)\v (x\w z)\\
x\v (y\w z)&=(x\v y)\w (x\v z).
\end{align}
\end{defn}

\begin{defn}
An element $z$ of a lattice $\L$ is called\/ {\em
join-irreducible} (respectively\/ {\em meet-irreducible}) if
$z=x\v y$ (respectively $z=x\w y$) implies $z=x$ or $z=y$. The set
of join-irreducible (respectively meet-irreducible) elements of
$\L$ is denoted by $J\left(\L\right)$ (respectively
$M\left(\L\right)$), or just by $J$ (respectively $M$) if
 no confusion is possible.
\end{defn}

\begin{defn}\label{irred}
An element in $J\left(\L\right)\cap M\left(\L\right)$ is called
\emph{irreducible}.

\medskip In the sequel, we shall denote  $J\left(\L\right)\cap
M\left(\L\right)$ by $JM\left(\L\right)$, or just $JM$ if no
confusion is possible.
\end{defn}
\begin{defn}
A subset $I$ of a poset $P$ is called an {\em ideal} of $P$ if for
all $x,\,y\in P$,
$$x\in I\text{ and }y\le x\text{ imply }y\in I.$$
\end{defn}

\begin{thm}[Birkhoff]\label{5.10}
Let $\L$ be a distributive lattice with $\hat{0}$, and $P$ the
poset of its nonzero join-irreducible elements. Then $\L$ is
isomorphic to the lattice of ideals of $P$, by means of the
lattice isomorphism
$$\a\mapsto I_{\a}:=\{\t\in P\mid \t\le\a\},\qquad \a\in\L.$$
\end{thm}

The following Lemma is easily checked.
\begin{lem}
With the notations as above, we have

$(a)$ $J=\{\t\in\L\mid\text{there exists at most one cover of the
form } (\t,\l)\}$.

$(b)$ $M=\{\t\in\L\mid\text{there exists at most one cover of the
form } (\l,\t)\}$.
\end{lem}

\begin{lem}[cf. \cite{LM}]\label{cover}
Let $(\t,\l)$ be a cover in ${\mathcal{L}}$. Then $I_{\t}$ equals
$I_{\l}\dot\cup\{\beta\}$ for some $\beta\in
J\left({\mathcal{L}}\right)$.
\end{lem}

\subsection{The variety $X(\L)$}\label{s3}

Consider the polynomial algebra $K[X_\alpha,\alpha\in\L]$; let
$\mathfrak{a}(\L)$ be the ideal generated by $\{X_\alpha
X_\beta-X_{\alpha\v \beta}X_{\alpha\w \beta},
\alpha,\beta\in\L\}$. Then one knows (cf.\cite{Hi}) that
$K[X_\alpha,\alpha\in\L]\,/\mathfrak{a}(\L)$ is a normal domain; in
particular, we have that $\mathfrak{a}(\L)$ is a  prime ideal. Let
$X(\L)$ be the affine variety of the zeroes in $K^l$ of
$\mathfrak{a}(\L)$ (here, $l=\# \L$).
 Then $X(\L)$ is an affine normal variety defined by binomials. On the other hand,
by \cite{ES}, we have that a binomial prime ideal is toric (here,
``toric ideal" is in the sense of \cite{St}, Chapter 4). Hence
$X(\L)$ is a toric variety for the action by a suitable torus $T$.

In the sequel, we shall denote
$R\left(\mathcal{L}\right):=K[X_\alpha,\alpha\in\L]\,/\mathfrak{a}(\L)$.
Further, for $\alpha\in \L$, we shall denote the image of
$X_\alpha$ in $R\left(\mathcal{L}\right)$ by $x_{\alpha}$.

\begin{defn} The variety $X\left(\mathcal{L}\right)$ will be called
a \em{Hibi toric variety.}
\end{defn}

\begin{remark} An extensive study of $X\left(\mathcal{L}\right)$ appears first in \cite{Hi}.
\end{remark}

We have that $\dim X(\L)=\dim T$.

\begin{thm}[cf. \cite{LM}]\label{dim}
The dimension of $X(\L)$ is equal to  $\#J\left(\L\right)$.
Further, dim$\,X\left({\L}\right)$ equals the cardinality of the
set of elements in a maximal chain in (the graded poset) ${\L}$.
\end{thm}

\subsection{Cone and dual cone of $X\left({\mathcal{L}}\right)$}
\label{cone}

As above, denote the poset of join-irreducibles in ${\L}$ by
$J\left({\L}\right)$ or just $J$. Denote by
${\mathcal{I}}\left({J}\right)$ the poset of ideals of $J$. For
$A\in {\mathcal{I}}\left({J}\right)$, denote by ${\mathbf{m}}_A$
the monomial:
$${\mathbf{m}}_A:={\underset{\tau\in A}{\prod}} \,y_\tau$$ in the
polynomial algebra $K[y_\tau,\tau\in J\left({\L}\right)]$. If
$\alpha$ is the element of $\L$ such that $I_\alpha=A$ (cf.
Theorem \ref{5.10}), then we shall denote ${\mathbf{m}}_A$ also by
${\mathbf{m}}_\alpha$. Consider the surjective algebra map
$$F:K[X_\alpha,\alpha\in{\mathcal{L}}]\rightarrow
K[{\mathbf{m}}_A,A\in{\mathcal{I}}\left({J}\right)],\,X_\alpha\mapsto
{\mathbf{m}}_A,\ A=I_\alpha
$$

\begin{thm}[cf. \cite{Hi}, \cite{LM}]\label{main1} We have an isomorphism
$$K[X\left({\L}\right)]\cong
K[{\mathbf{m}}_A,A\in{\mathcal{I}}\left({J}\right)].$$
\end{thm}

Let us denote the torus acting on the toric variety
$X\left({\L}\right)$ by $T$; by Theorem \ref{dim}, we have,
$\dim\,T=\# J\left({\L}\right)=d$, say. Identifying $T$ with
$(K^*)^d$, let $\{f_z,z\in J\left({\L}\right)\}$ denote the
standard ${\mathbb{Z}}$-basis  for $X(T)$, namely, for $t=(t_z, z
\in J\left({\L}\right))$, $f_z(t)=t_z$. Denote $M:=X(T)$; let $N$
be the ${\mathbb{Z}}$-dual of $M$, and $\{e_y,y\in
J\left({\L}\right)\}$ be the basis of $N$ dual to $\{f_z,z\in
J\left({\L}\right)\}$. For $A\in {\mathcal{I}}\left({J}\right)$,
set
$$f_A:={\underset{z\in A}{\sum}} \,f_z$$ Let
$V=N_{{\mathbb{R}}}(=N\otimes_{{\mathbb{Z}}} {\mathbb{R}})$. Let
$\sigma \subset V$ be the cone such that
$X\left({\mathcal{L}}\right)=X_\sigma$.

As an immediate consequence of Theorem \ref{main1}, we have
\begin{prop}\label{semi}
The semigroup $S_\sigma$ is generated by $f_A, A\in
{\mathcal{I}}\left({J}\right)$.
\end{prop}

Let $M(J\left({\L}\right))$ be the set of maximal elements in the
poset $J\left({\L}\right)$. Let $Z(J\left({\L}\right))$ denote the
set of all covers in the poset $J\left({\L}\right)$. For a cover
$(y,y')\in Z(J\left({\L}\right))$, denote
$$v_{y,y'}:=e_{y'}-e_{y}$$

\begin{prop}[cf. \cite{LM}, Proposition 4.7]\label{prop-cpc} The cone $\sigma$ is
generated by $\{e_{z},z\in M(J\left({\L}\right)),\,v_{y,y'},
(y,y')\in Z(J\left({\L}\right))\}$.
\end{prop}

\subsection{The sublattice $D_\tau$}\label{anal}
We shall concern ourselves
 just with the closed points in $X\left( {\L}\right)$.
So in the sequel, by a point in $X\left( {\L}\right)$, we shall
mean a closed point. Let $\tau$ be a face of $\sigma$, and
$P_{\tau}$ the distinguished point (cf. \S \ref{distinguished})

For a point $P\in X\left( {\L}\right)$ (identified with a point in
${\mathbb{A}}^l$, $l=\#\L$), let us denote by $P(\alpha)$, the
$\alpha$-th co-ordinate of $P$. Let
$$D_{\tau}=\{\alpha\in{\L}\,|\,P_{\tau}(\alpha)\neq 0\}$$ We have,
\begin{lem}[cf. \cite{LM}]
$D_{\tau}$ is an embedded sublattice.
\end{lem}

Conversely, we have
\begin{lem}[cf. \cite{LM}]\label{embed}
Let $D$ be an embedded sublattice in ${\L}$. Then $D$ determines a
unique face $\tau$ of $\sigma$ such that $D_{\tau}$ equals $D$.
\end{lem}

Thus in view of the two Lemmas above, we have a bijection
$$\{\mathrm{\ faces\ of\
}\sigma\}\,{\buildrel{bij}\over{\leftrightarrow}}\, \{\mathrm{\
embedded\ sublattices \ of\ }{\L}\}$$
\begin{prop}[cf. \cite{LM}]
Let $\tau$ be a face of $\sigma$. Then we have
${\overline{O_{\tau}}}=X\left({D_{\tau}}\right)$.
\end{prop}


\section{Grid Lattices}\label{grid}

In this section, we restrict our attention to a specific class of
distributive lattices, and show that some desirable properties
hold.  Give $\N\times\N$ the lattice structure
\[ (\alpha_1,\alpha_2 )\wedge (\beta_1,\beta_2)=(\delta_1,\delta_2),\ (\alpha_1,\alpha_2 )\vee (\beta_1,\beta_2 )=(\gamma_1,\gamma_2),\]
where $\delta_i= $ min$\{\alpha_i,\beta_i\}$, $\gamma_i=$
max$\{\alpha_i,\beta_i\}$.

\begin{definition} Let $J$ be a finite, distributive sublattice of $\N\times\N$, such that if $\alpha$ covers $\beta$ in $J$, then $\alpha$ covers $\beta$ in $\N\times\N$ as well.  Then we say $J$ is a \emph{grid lattice}.
\end{definition}

\begin{remark}\label{remark-6-2} For $J$ a grid lattice, we have the following:
\begin{enumerate}
\item $J$ is a distributive lattice. \item For any $\mu\in J$,
there exist at most two distinct covers of the form $(\alpha,\mu
)$ in $J$, i.e., there are at most two elements in $J$ covering
$\mu$. \item For any $\lambda\in J$, $\lambda$ covers at most two
distinct elements in $J$. \item If $\alpha$, $\beta$ are two
covers of $\mu$ in $J$, then $\alpha\vee\beta$ covers both
$\alpha$, $\beta$; thus the interval $[\mu , \alpha\vee\beta ]$ is
a rank 2 subposet of $J$.
\end{enumerate}
\end{remark}

\noindent\textbf{Example:}
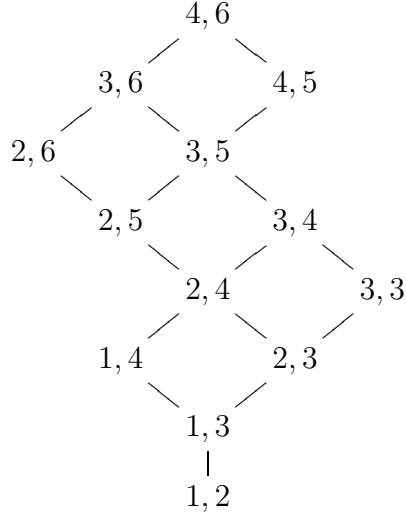
\[ \xymatrix@-15pt{
& & 4,6\ar@{-}[dr]\ar@{-}[dl]& & &  \\
& 3,6 \ar@{-}[dr]\ar@{-}[dl]& & 4,5\ar@{-}[dl] & & \\
2,6 \ar@{-}[dr]& & 3,5\ar@{-}[dr]\ar@{-}[dl] & & & \\
& 2,5\ar@{-}[dr]& & 3,4\ar@{-}[dr]\ar@{-}[dl] & &  \\
& & 2,4 \ar@{-}[dr]\ar@{-}[dl]& & 3,3\ar@{-}[dl] &  \\
& 1,4\ar@{-}[dr]& & 2,3\ar@{-}[dl] & &  \\
& & 1,3\ar@{-}[d]& & & \\
& & 1,2& & & }\]

\subsection{} For the rest of this section, let $J$ be a grid
lattice, and let $\mathcal{L}$ be the poset of ideals of $J$. From
Theorem \ref{5.10}, we have that $\mathcal{L}$ is a distributive
lattice with $J$ as its poset of join irreducibles. Thus we will
correlate join irreducible elements in $\L$ with elements of $J$.
Recall that for $x,y\in\L$, $x\geq y$ if and only if $I_x\supseteq
I_y$ as ideals in $J$.

\begin{lemma}\label{lemma-a1} Given $\gamma_1,\gamma_2\in J$, $(\gamma_1\wedge\gamma_2 )_{\L}$ belongs to $J$ and is in fact equal to $(\gamma_1\wedge\gamma_2 )_J$.
\end{lemma}
\begin{proof} Let $\theta =(\gamma_1\wedge\gamma_2 )_J$ and
$\phi =(\gamma_1\wedge\gamma_2)_{\L}$.  Clearly $\theta\in
I_{\gamma_1}\cap I_{\gamma_2}=I_\phi$.  Therefore $I_\theta\subset
I_\phi$.  Let now $\eta\in I_\phi(\subset J)$.  Then
$\eta\leq\phi$, and thus $\eta$ is less than or equal to both
$\gamma_1$ and $\gamma_2$ in $\L$, and therefore in $J$.  Hence
$\eta\leq\theta$, and thus $I_\phi\subset I_\theta$.  The result
follows.
\end{proof}

\begin{lemma}\label{lemma-a2} Let ($\alpha$, $\beta$) be an
incomparable pair of irreducibles (cf. Definition \ref{irred}) in
$\L$. Then
\begin{enumerate}
\item $\alpha$, $\beta$ are meet irreducibles in $J$, \item
$(\alpha\wedge\beta )_{\L}=(\alpha\wedge\beta )_J\in J$.
\end{enumerate}
\end{lemma}
\begin{proof} Part (2) follows from Lemma \ref{lemma-a1},
(note that $\alpha$, $\beta\in J$).  Now say $\alpha =
(\gamma_1\wedge \gamma_2 )_J$ for an incomparable pair
($\gamma_1$, $\gamma_2$) in $J$.  Lemma \ref{lemma-a1} implies
that $\alpha = (\gamma_1\wedge\gamma_2 )_{\L}$, a contradiction
since $\alpha$ is meet irreducible in $\L$.  Part (1) follows.
\end{proof}

Thus an incomparable pair $(\alpha,\beta)$ of irreducibles in
${\L}$ determines a (unique) non-meet irreducible in $J$ (namely,
$(\alpha\wedge\beta )_{\L}=(\alpha\wedge\beta )_{J}$). We shall
now show (cf. Lemma \ref{lemma-a3} below) that conversely a
non-meet irreducible element $\mu$ in $J$ determines a unique
incomparable pair $(\alpha,\beta)$ of irreducibles in ${\L}$. We
first prove a couple of preliminary results:
\begin{lem}\label{lem1}
Let $\mu$ be a non-meet irreducible element  in $J$. Then $\mu$
determines an incomparable pair $(\alpha,\beta)$ of elements (in
$J$) both of which are meet irreducible in $J$.
\end{lem}
\begin{proof}
Let $\mu=(\mu_1,\mu_2)$ (considered as an element of
$\N\times\N$). Since $\mu$ is non-meet irreducible element  in
$J$, there exist $x=(x_1,x_2), y=(y_1,y_2)$ in $J$, $x,y>\mu$ such
that $x_2>\mu_2, y_1>\mu_1$. Define
$\alpha=(\alpha_1,\alpha_2),\beta=(\beta_1,\beta_2)$ in $J$ as
\[\alpha = \mbox{ the maximal element $x>\mu$ in }J
\mbox{ such that }x_1=\mu_1 ,\]
\[\beta = \mbox{ the maximal element $y>\mu$ in }J\mbox{ such that }
y_2=\mu_2 .\] Clearly $\alpha,\beta$ are both meet-irreducible in
$J$ (note that $(\mu_1+1,\alpha_2)$ (resp. $(\beta_1,\mu_2+1)$) is
the unique element in $J$ covering $\alpha$ (resp. $\beta$) in
$J$). Also, it is clear that $(\alpha,\beta)$ is an incomparable
pair.
\end{proof}
Let $\mu,\alpha,\beta$ be as in the above Lemma. In particular, we
have, $\mu_1=\alpha_1<\beta_1,\,\mu_2=\beta_2<\alpha_2$.

\begin{lem}\label{maximal-remark} With notation as in
Lemma \ref{lem1}, we have,
\begin{enumerate}
\item $(\alpha\vee\beta)_J =( \beta_1,\alpha_2 )$. \item $\alpha$
is the maximal element of the set $\{x=(x_1,x_2 )\in J\mid
x_1=\alpha_1\}$, and $\beta$ is the unique maximal element of the
set $\{x=(x_1,x_2 )\in J\mid x_2=\beta_2\}$.
\end{enumerate}
\end{lem}
\begin{proof}
Assertion (2) is immediate from the definition of $\alpha,\beta$.
Assertion (1) is also clear.
\end{proof}
\begin{lemma}\label{lemma-a3} Let $\mu,\alpha , \beta$ be as in
Lemma \ref{lem1}. Then $\alpha$ and $\beta$ are irreducibles in
$\L$. Thus the non-meet irreducible element $\mu$ of $J$
determines a unique incomparable pair of irreducibles in $\L$.
\end{lemma}
\begin{proof} We will show the result for $\alpha$ (the proof for $\beta$ being similar).  Since $\alpha\in
J,\,\alpha$ is join irreducible in $\L$. It remains to show that
$\alpha$ is meet irreducible in $\L$.  If possible, let us assume
that there exists an incomparable pair $(\theta_1 ,$ $\theta_2)$ in $\L$
such that $\theta_1\wedge\theta_2=\alpha$; without loss of
generality, we may suppose that $\theta_1$ and $\theta_2$ both
cover $\alpha$. Then there exist (cf. Lemma \ref{cover})
$\gamma,\delta\in J$ such that
\[ I_{\theta_1}= I_\alpha\dot\cup\{\gamma\},\,I_{\theta_2}=
I_\alpha\dot\cup\{\delta\}.\] We have \[I_{\gamma}\cap
I_{\delta}\subset I_{\theta_1}\cap I_{\theta_2} =
I_\alpha . \eqno{(\ast )}\]   Also, $\gamma$, $\delta$ are either covers
of $\alpha$ in $J$, or non-comparable to $\alpha$. (They cannot be
less than $\alpha$ because they are not in $I_\alpha$.)

\textbf{Case 1:} Suppose $\gamma$ and $\delta$ are covers of
$\alpha$ in $J$.  Then $\alpha$ is not meet irreducible in $J$, a
contradiction (cf. Lemma \ref{lemma-a2},(1)).

\textbf{Case 2:} Suppose $\gamma$ covers $\alpha$ in $J$, and
$\delta$ is non-comparable to $\alpha$. Let $\delta =
(\delta_1,\delta_2 ),\,\xi = (\xi_1,\xi_2 )=(\alpha\vee\delta
)_J$. Then the fact that $\xi>\alpha$ (since $\alpha,\delta$ are
incomparable) implies (in view of Lemma \ref{maximal-remark}, (2))
that $\xi_1>\mu_1$; hence $\delta_1(=\xi_1)\ge \mu_1+1$, and
$\delta_2<\alpha_2$. Also, $\gamma =(\mu_1+1,\alpha_2)$ (cf. Lemma
\ref{maximal-remark}, (2)). Therefore $\gamma\wedge\delta =
(\mu_1+1,\delta_2)$, but this element is non-comparable to
$\alpha$, and thus $I_\gamma\cap I_\delta\not\subset I_\alpha$, a
contradiction to $(\ast )$. Hence we obtain that the possibility
``$\gamma$ covers $\alpha$ in $J$ and $\delta$ is non-comparable
to $\alpha$" does not exist.  A similar proof shows that the possibility  
``$\delta$ covers $\alpha$ in $J$ and $\gamma$ is non-comparable to $\alpha$" does
not exist.

\textbf{Case 3:} Suppose both $\gamma = (\gamma_1,\gamma_2)$ and
$\delta=(\delta_1,\delta_2)$ are non-comparable to
$\alpha=(\mu_1,\alpha_2)$.  As in Case 2, we must have $\delta_2
<\alpha_2$, and thus $\delta_1 >\mu_1$. Similarly, $\gamma_2
<\alpha_2$, $\gamma_1 >\mu_1$. Thus the minimum of
$\{\gamma_1,\delta_1\}$ is still greater than $\mu_1$, therefore
$I_\gamma\cap I_\delta\not\subset I_\alpha$,
 a contradiction to $(\ast )$.

 Thus our assumption that $\alpha$ is non-meet irreducible in $\L$
 is wrong, and it follows that $\alpha$ (and similarly $\beta$)
 is meet irreducible in $\L$.
\end{proof}

We continue with the above notation; in particular, we denote
$\mu=(\mu_1,\mu_2),\mu_1=\alpha_1<\beta_1,\,\mu_2=\beta_2<\alpha_2$.

\begin{lemma}\label{V-shape} Let $x=(x_1,x_2)\in J$.
 If $x\not\in I_\alpha\cup I_\beta$, then $x >
\alpha\wedge\beta$.
\end{lemma}
\begin{proof} By hypothesis, we have $x\not\leq\alpha$,
$x\not\leq\beta$.

We first claim that $x_1
>\alpha_1$; for, if possible, let us assume $x_1\leq\alpha_1$.
Since $x\not\leq\alpha$, we must have $x_2
>\alpha_2$.  Thus $x\vee\alpha= (\alpha_1,x_2 )$
($>\alpha$, since $\alpha\not\ge x$); but this is a contradiction,
by the property of $\alpha$ (cf. Lemma \ref{maximal-remark},(2)).
Hence our assumption is wrong, and we get $x_1
>\alpha_1$.

Similarly, we have, $x_2>\beta_2$, and the result follows (note
that by our notation (and definition of $\alpha,\beta$), we have
$\alpha\wedge\beta= (\alpha_1,\beta_2)$).
\end{proof}

\begin{definition}\label{lab} For an incomparable (unordered) pair $(\alpha,\beta)$ of
irreducible elements in $\L$, define
 $$\L _{\alpha,\beta}=\L\setminus
[\alpha\wedge\beta,\alpha\vee\beta ].$$
\end{definition}

\begin{prop}\label{embedded} $\L_{\alpha,\beta}$ is an embedded sublattice.
\end{prop}
\begin{proof}  First, we show that $\L_{\alpha,\beta}$ is a
sublattice.  To do this, we identify $\L$ with the ``lattice of
ideals'' of $J $.  Thus, for $x\in \L_{\alpha,\beta}$, either
$I_x\not\supset (I_\alpha\cap I_\beta)$ or $I_x\not\subset
(I_\alpha\cup I_\beta )$, by definition of $\L _{\alpha,\beta}$.
Note that $I_\alpha\cap I_\beta=I_{\alpha\wedge\beta}$, and
$I_\alpha \cup I_\beta = I_{\alpha\vee\beta}$.

\textbf{Case 1:} Let $x,y\in\L_{\alpha, \beta}$ such that
$I_x,\,I_y\not\supset I_{\alpha\wedge\beta}$.  Then clearly
$I_x\cap I_y\not\supset I_{\alpha\wedge\beta}$; and thus $x\wedge
y\in \L_{\alpha,\beta}$.  We also have (by the definition of
ideals) that $\alpha\wedge\beta\not\in I_x,I_y$ (note that
$\alpha\wedge\beta\in J$ (cf. Lemma \ref{lemma-a2},(2))),
therefore $\alpha\wedge\beta\not\in I_x\cup I_y$, and therefore
$x\vee y\in \L_{\alpha,\beta}$.

\textbf{Case 2:} Let $x,y\in\L_{\alpha,\beta}$ such that
$I_x\not\supset I_{\alpha\wedge\beta}$ and $I_y\not\subset
I_{\alpha\vee\beta}$.  Then clearly $I_x\cap I_y\not\supset
I_{\alpha\wedge\beta}$ and $I_x\cup I_y\not\subset I_{\alpha}\cup
I_\beta$. Hence, $x\vee y,x\wedge y$ are in $\L_{\alpha,\beta}$.

\textbf{Case 3:} Let $x,y\in \L_{\alpha,\beta}$ such that $I_x,\,
I_y\not\subset I_{\alpha\vee\beta}$.
 Clearly $I_x\cup I_y\not\subset I_{\alpha}\cup I_\beta$; hence,
 $x\vee y\in\L_{\alpha,\beta}$.

 \noindent\textbf{Claim:} $I_x\cap I_y\not\subset I_\alpha\cup
 I_\beta$. 
 
 Note that Claim implies that $x\wedge y\in\L_{\alpha,\beta}$.
  If possible, let us assume that
 $I_x\cap I_y\subset I_\alpha\cup I_\beta$. Now
 the hypothesis that $I_x,\, I_y\not\subset I_{\alpha}\cup
 I_\beta$ implies that
 there exist $\theta,\delta\in J$ such that $\theta\in I_x,\theta\not\in I_\alpha\cup I_\beta$,
  and $\delta\in I_y,\delta\not\in I_\alpha\cup I_\beta$. Now
  $I_\theta\cap I_\delta\subset I_x\cap I_y\subset I_\alpha\cup
   I_\beta$ (note that by our assumption, $I_x\cap I_y\subset I_\alpha\cup
   I_\beta$). Hence we obtain that either
    $\theta\wedge\delta \leq \alpha$ or
    $\theta\wedge\delta\leq\beta$; let us suppose $\theta\wedge\delta
    \le \alpha$ (proof is similar if $\theta\wedge\delta\leq\beta$).  By Lemma \ref{V-shape}, we have that both
     $\theta ,\delta\geq \alpha\wedge\beta$, and
     hence $\theta\wedge\delta\geq\alpha\wedge\beta$.  Thus
\[\alpha\geq\theta\wedge\delta\geq\alpha\wedge\beta=
(\alpha_1,\beta_2 ) \eqno{(\ast\ast )}\] Let
$\xi(=(\xi_1,\xi_2))=\theta\wedge\delta $.  Then $(\ast\ast)$ implies that
$\xi_1 = \alpha_1$; hence at least one of $\{\theta_1,\delta_1\}$,
say $\theta_1$ equals $\alpha_1$. This implies that $\theta_2
>\alpha_2$ (since $\theta\not\in I_\alpha$).  This contradicts
Lemma \ref{maximal-remark}(2). Hence our assumption is wrong and
it follows that $I_x\cap I_y\not\subset I_\alpha\cup
 I_\beta$.

This completes the proof in  Case 3. Thus we have shown that
$\L_{\alpha,\beta}$ is a sublattice.

Next, we will show that $\L_{\alpha,\beta}$ is an embedded
sublattice.  Let $x,y\in \L$ be such that $x\vee y,x\wedge y$ are
in $\L_{\alpha,\beta}$. We need to show that $x,y \in
\L_{\alpha,\beta}$. This is clear if either $x\wedge
y\not\le\alpha\vee\beta$ or $x\vee y\not\ge\alpha\wedge\beta$ (in
the former case, $x,y\not\le\alpha\vee\beta$, and in the latter
case, $x,y\not\ge\alpha\wedge\beta$). Let us then suppose that
$x\wedge y\le\alpha\vee\beta$ and $x\vee y \ge\alpha\wedge\beta$;
this implies that $x\wedge y\not\ge\alpha\wedge\beta$ and $x\vee
y\not\le\alpha\vee\beta$ (since, $x\vee y,x\wedge y$ are in
$\L_{\alpha,\beta}$), i.e., $I_x\cap I_y\not\supset I_{\alpha}\cap
I_\beta$ and $I_x\cup I_y\not\subset I_\alpha\cup I_\beta$.  We
will now show that $x,y\in\L_{\alpha,\beta}$.

Since $\alpha\wedge\beta\not\in I_x\cap I_y$, we have that one of
the elements $\{x,y\}$ must not be greater than or equal to
$\alpha\wedge \beta$, say $x\not\geq\alpha\wedge\beta$.
 This implies that $x\in\L_{\alpha,\beta}$.  It remains to
 show that $y\in\L_{\alpha,\beta}$.
 If $y\not\geq \alpha\wedge\beta$, then we would obtain that
 $y\in\L_{\alpha,\beta}$.  Let us then assume that $y\geq\alpha\wedge\beta$;
  i.e. $I_y\supset I_\alpha\cap I_\beta$.  Note that for any
  $\delta\in I_x$, we have $\delta\leq x$ and thus
  $\delta\not\geq\alpha\wedge\beta$.  By Lemma \ref{V-shape},
  $\delta\in I_\alpha\cup I_\beta$, and
  therefore $I_x\subset I_\alpha\cup I_\beta$.  Since by hypothesis
   $I_x\cup I_y\not\subset I_\alpha\cup I_\beta$,
   we must have $I_y\not\subset I_\alpha\cup I_\beta$.
   Therefore, $y\in \L_{\alpha,\beta}$.

   This completes the proof of the assertion that $\L_{\alpha,\beta}$ is an embedded
sublattice, and therefore the proof of the Proposition.
\end{proof}

\section{Singular locus of $X({\L})$}\label{faces} In this section, we
determine the singular locus of $X({\L})$, ${\L}$ being as in \S
\ref{grid}. Let $\sigma$ be the cone associated to
$X\left(\L\right)$. We follow the notation of \S \ref{first-sect}
and \S \ref{varlat}.

\begin{definition}\label{def-sing-face} A face $\tau$ of $\sigma$ is a \emph{singular} (resp. \emph{non-singular}) face if $P_\tau$ is a singular (resp. non-singular) point of $X_{\sigma}$.
\end{definition}

\begin{definition}\label{gen} Let us denote by $W$ the set of generators for $\sigma$ as
described in Proposition \ref{prop-cpc}.  Let $\tau$ be a face of
$\sigma$, and let $D_\tau$ be as in \S\,\ref{anal}.  Define
\[W\left(\tau\right)=\{v\in W\mid f_{I_{\alpha}}
\left(v\right)=0,\,\forall\,\alpha\in D_\tau\}.\] Then
$W\left(\tau\right)$ gives a set of generators for $\tau$.
\end{definition}

\subsection{Determination of $W\left(\tau\right)$}\label{Wij}
Let $(\alpha,\beta)$ be an incomparable (unordered) pair of
irreducible elements of ${\L}$. By Proposition \ref{embedded},
${\L} _{\alpha,\beta}$ is an embedded sublattice of ${\L}$ (${\L}
_{\alpha,\beta}$ being as in Definition \ref{lab}).  Let
$\tau_{\alpha ,\beta}$ be the face of $\sigma$ corresponding to
${\L} _{\alpha,\beta}$ (cf. Lemma \ref{embed}; note that
$D_{\tau_{\alpha, \beta}}=\L _{\alpha,\beta}$). Let us denote
$\tau=\tau_{\alpha,\beta}$. Following the notation of \S
\ref{grid}, let
$\mu(=(\mu_1,\mu_2))=\alpha\wedge\beta,\alpha_1=\mu_1,\beta_2=\mu_2$.
 Since $\mu$ is not meet irreducible in $J$, there are two elements
  $A$ and $B$ in $J$ covering $\mu$, namely,
$A=(\alpha_1,\beta_2+1),B=(\alpha_1+1,\beta_2)$. Also, we have
that $A\vee B$ (in the lattice $J$)  covers both $A$ and $B$, (cf.
Remark \ref{remark-6-2}). Let $C=(A\vee B)_J$; then
$C=(\alpha_1+1,\beta_2+1)$.

It will aid our proof below to notice a few facts about
 the generating set  $W(\tau)$ of $\tau$.  First of all, $e_{\hat{1}}$ is not a
generator for any $\tau_{\alpha,\beta}$; because $\hat{1}\in \L
_{\alpha,\beta}$ for all pairs $(\alpha,\beta )$, and
$e_{\hat{1}}$ is non-zero on $f_{I_{\hat{1}}}$.

Secondly, for any cover $(y,y'),y>y'$ in $J\left({\L}\right)$,
$e_{y'}-e_{y}$ is not a generator of $\tau$ if
$y'\in\mathcal{L}_{\alpha,\beta}$, because
$f_{I_{y'}}\left(e_{y'}-e_{y}\right)\neq 0$.  Thus, in
determining the elements of $W\left(\tau\right)$, we need only be
concerned with elements $e_{y'}-e_{y}$ of $W$ such that $y'\in
J\cap [\alpha\wedge\beta,\alpha\vee\beta]$.

\begin{lemma}\label{lemma-a4} $J\cap [\alpha\wedge\beta,\alpha\vee\beta]=\{x\in J\mid x\in [\mu,\alpha ]\cup [\mu,\beta]\}$.
\end{lemma}
\begin{proof} The inclusion $\supseteq$ is clear.
To show the inclusion $\subseteq$, let $x\in J\cap
[\alpha\wedge\beta,\alpha\vee\beta]$. If possible, assume
$x\not\in [\mu,\alpha]\cup [\mu,\beta ]$; the assumption implies
that $x\not\in I_\alpha\cup I_\beta(=I_{\alpha\vee\beta})$. Hence
we obtain that $x\not\le\alpha\vee\beta$, a contradiction to the
hypothesis that $x\in [\alpha\wedge\beta,\alpha\vee\beta]$.
\end{proof}

\begin{lemma}\label{unique-C} The set $\{x\in J\mid x\not\in I_\alpha\cup I_\beta\}$ has a unique minimal element; moreover that element is $C$.
\end{lemma}
\begin{proof}  For any $x$ in this set, we have
$x> \alpha\wedge\beta$ (cf. Lemma \ref{V-shape}). Hence by Lemma
\ref{maximal-remark},(2), and the hypothesis that $x\not\in
I_\alpha\cup I_\beta$, we obtain that $x_1>\alpha_1$,
$x_2>\beta_2$. Therefore,
\[ \{x\in J\mid x\not\in I_\alpha\cup I_\beta\}=\{x\in J\mid x_1 >\alpha_1,\, x_2 >\beta_2\}.\]
This set clearly has a minimal element, namely
$C=(\alpha_1+1,\beta_2+1)$.
\end{proof}

\begin{theorem}\label{main} Following the notation from above, we
have
\[W\left(\tau\right)=\{e_\mu-e_A,\, e_\mu-e_B,\, e_A -e_C,\, e_B-e_C\}.\]
\end{theorem}
\begin{proof}

\noindent\textbf{Claim 1:}  $W\left(\tau\right)\supset
\{e_\mu-e_A,\, e_\mu-e_B,\, e_A -e_C,\, e_B-e_C\}$.

We must show that for any $x\in \L_{\alpha,\beta}$, $f_{I_x}$ is
zero on these four elements of $W$.  If possible, let us assume
that there exists a $x\in \L_{\alpha,\beta}$  such that $f_{I_x}$
is non-zero on some of the above four elements.  Then clearly
$x\geq \mu (= \alpha\wedge \beta)$. Hence
$x\not\leq\alpha\vee\beta$ (since $x\not\in
[\alpha\wedge\beta,\alpha\vee\beta]$), i.e., $I_x\not\subset
I_\alpha\cup I_\beta$.  Therefore $I_x$ contains some join
irreducible $\gamma$ such that $\gamma\not\leq \alpha,\,\beta$;
hence, $I_\gamma\not\subset I_\alpha\cup I_\beta$. This implies
(cf. Lemma \ref{unique-C}) that $\gamma\geq C$. Hence we obtain
that $C\in I_x$. Therefore, $x\geq C$, and $f_{I_x}$ is zero on
all of the four elements of Claim 1, a contradiction to our
assumption. Hence our assumption is wrong and Claim 1 follows.

\noindent\textbf{Claim 2:}  $W\left(\tau\right)=\{e_\mu-e_A,\,
e_\mu-e_B,\, e_A -e_C,\, e_B-e_C\}$.

In view of $\S$ \ref{Wij}, it is enough to show that for all
$\theta\in J\cap [\alpha\wedge\beta,\alpha\vee\beta]$, the element
$e_{\theta}-e_{\delta}\in W$ which is different from the four
elements of Claim 1 is not in $W(\tau )$. In view of Lemma
\ref{lemma-a4}, it suffices to examine all covers in $J$ of all
elements in $ \left( [\mu,\alpha ]\cup [\mu,\beta ]\right)_J$.
This diagram represents the part of the grid lattice $J$ we are
concerned with:
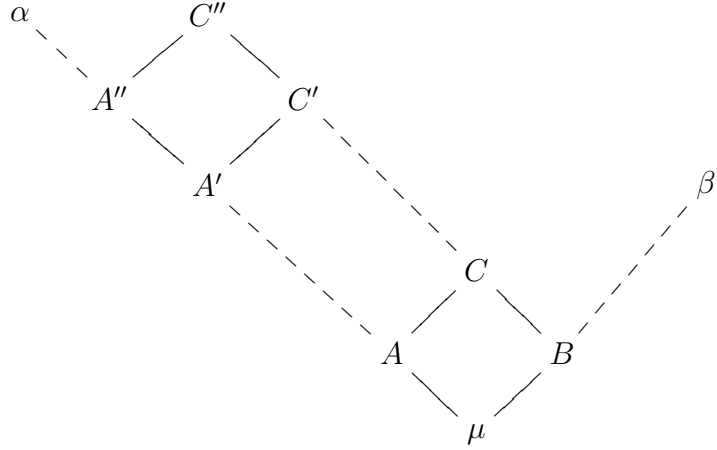
\[ \xymatrix@-8pt{
\alpha & & C''& & & & & &\\
& A'' \ar@{-}[ur]\ar@{--}[ul]& & C'\ar@{-}[ul] & & & & &\\
& & A'\ar@{-}[ur]\ar@{-}[ul] & & & & & &\beta \\
& & & & & C\ar@{--}[uull] & & &\\
& & & & A\ar@{--}[uull]\ar@{-}[ur] & & B\ar@{-}[ul]\ar@{--}[uurr]& & \\
& & & & & \mu \ar@{-}[ur]\ar@{-}[ul]& & &
 } \]

In the diagram above, consider $A'=(\alpha_1,\beta_2 +n)$,
$A''=(\alpha_1,\beta_2+n+1)$, and $C'=(\alpha_1+1,\beta_2+n)$.
Note that all elements of $J\cap [\mu,\alpha ]$ can be written in
the form of $A'$.  Thus, we need to check elements
$e_{A'}-e_{A''}$ and $e_{A'}-e_{C'}$ in $W$.

First, we observe that $C'\in\L_{\alpha,\beta}$, and $f_{I_{C'}}$
is non-zero on $e_{A'}-e_{A''}$. Hence $e_{A'}-e_{A''}\not\in
W(\tau)$.

 Next, let $x= (A'\vee C)_{\L}$, (note that $x$ is not in
$J$, and thus does not appear on the diagram above). Then
$I_x=I_{A'}\cup I_C$; and we have $x\in \L_{\alpha,\beta}$ (since,
$C\not\in I_\alpha\cup I_\beta$ and $x>C$, we have,
$x\not\le\alpha\vee\beta$). Moreover, $f_{I_x}$ is non-zero on
$e_{A'}-e_{C'}$. Hence $e_{A'}-e_{C'}\not\in W(\tau)$.

This completes the proof for the interval $[\mu, \alpha]$, and a
similar discussion yields the same result for the interval $[\mu
,\beta ]$.

Thus Claim 2 (and hence the Theorem) follows.
\end{proof}
As an immediate consequence of Theorem \ref{main}, we have the
following
\begin{theorem}\label{quadric} Let $(\alpha,\beta)$ be an incomparable pair of
irreducibles in ${\L}$. We have an identification of the (open)
affine piece in $X(\L)$ corresponding to the face
$\tau_{\alpha,\beta}$ with the product $Z\times (K^*)^{\#{\L}-3}$,
where $Z$ is the cone over the quadric surface $x_1x_4-x_2x_3=0$
in $\mathbb{P}^3$.
\end{theorem}

\begin{lemma}\label{dim-tau-ab} The dimension of the face
$\tau_{\alpha,\beta}$ equals 3.\end{lemma}
\begin{proof}  By Theorem \ref{main}, a set of generators for
$\tau_{\alpha,\beta}$ is given by

\noindent $\{e_{\mu}-e_A,\,e_{\mu}-e_B,\,e_A-e_C,\,e_B-e_C\}$. We
see that a subset of three of these generators is linearly
independent. Thus if the fourth generator can be put in terms of
the first three, the result follows.  Notice that
\[\left(e_{\mu}-e_A\right)-\left(e_{\mu}-e_B\right)+\left(e_A-e_C\right)=e_B-e_C .\]
\end{proof}

The following fact from \cite{G-H-paper} (Lemmas 6.21, 6.22) holds for a general
toric variety.

\begin{fact}\label{basis-sing} Let $\tau$ be a face of $\sigma$.
Then $P_{\tau}$ is a smooth point of $X_{\sigma}$ if and only if
$X_\tau$ is non-singular.
\end{fact}

Combining the above fact with Theorem \ref{quadric}, we obtain
the following:
\begin{theorem}\label{sing}
$P_\tau\in $Sing$\, X_\sigma$, for $\tau =\tau_{\alpha,\beta}.$
Further, the singularity at $P_\tau$ is of the same type as that
at the vertex of the cone over the quadric surface
$x_1x_4-x_2x_3=0$ in $\mathbb{P}^3$.
\end{theorem}

Next, we will show that the faces containing some
$\tau_{\alpha,\beta}$ are the only singular faces.

\begin{lemma}\label{s-3} Let $(y,y'),y>y'$ be a cover in $J$.
Then either $e_{y'}-e_y\in W\left(\tau_{\alpha,\beta}\right)$ for
some incomparable pair $(\alpha,\beta)$ of irreducibles in ${\L}$,
or $y,y'$ are comparable to every other element of $J$.
\end{lemma}
\begin{proof} \textbf{Case 1:} Let $y'$ be non-meet
irreducible in $J$.

In view of the hypothesis, we can find an incomparable pair
$(\alpha,\beta)$ of irreducibles in ${\L}$ such that
$y'=\alpha\wedge\beta$, as shown in Lemma  \ref{lemma-a3} (with
$\mu=y'$). Thus $e_{y'}-e_y=e_{\mu}-e_{A}$ or $e_{y'}-e_y=e_\mu -
e_B$ as in Theorem \ref{main}.

\noindent\textbf{Case 2:} Let $y'$ be meet irreducible, but not
join irreducible (in $J$).

Let $x_1$ and $x_2$ be the two elements covered by $y'$ in $J$
(cf. Remark \ref{remark-6-2}); thus $(x_1\vee x_2)_J = y'$.

For convenience of notation, all join and meet operations in this
proof will refer to the join and meet operations in the lattice
$J$.

\textbf{Claim (a):} If both $x_1$ and $x_2$ are meet irreducible
(in $J$), then $y' , y$ are comparable to every element of $J$.

  If possible, let us assume that there exists a $z\in J$ such that $z$ is non-comparable to $y'$.
  We first observe that $z$ is non-comparable to both $x_1$ and
  $x_2$; for, say $z,x_1$ are comparable, then $z>x_1$ necessarily
  (since $z,y'$ are non-comparable). This implies that $x_1\le z\wedge y'
  <y'$, and hence we obtain that $x_1= z\wedge y'
  <y'$ (since $(y', x_1)$ is a cover),
  a contradiction to the hypothesis that $x_1$ is meet
  irreducible. Thus we obtain that $z$ is non-comparable to both $x_1$ and
  $x_2$.
  Now, we have, $z\vee x_i \geq z\vee y'$ (note that $x_i, i=1,2$ being meet
  irreducible in $J$, $y'$ is the unique element covering
  $x_i, i=1,2$, and hence $z\vee x_i \geq  y'$).
   Hence $(z\vee y'\ge)z\vee x_i\geq z\vee y'$,
   and we obtain
\[z\vee x_1 = z\vee y' = z\vee x_2.\]
On the other hand, the fact that $z\wedge y' < y'$ implies that
$z\wedge y' \le x_1$ or $x_2$.  Let $i$ be such that $z\wedge y'
\le x_i$. Then $z\wedge y' \leq z\wedge x_i \leq z\wedge y'$;
therefore
\[z\wedge x_i = z\wedge y' .\]
Now
\[ y'\wedge\left(x_i\vee z\right)=y'\wedge\left(y'\vee z\right)= y';\quad \left(y'\wedge x_i\right)\vee\left( y'\wedge z\right)=x_i\vee\left(x_i\wedge z\right)=x_i .\]
Therefore $J$ is not a distributive lattice (Definition
\ref{dist-def}), a contradiction. Hence our assumption is wrong
and it follows that $y'$ is comparable to every element of $J$,
and since $y$ is the unique cover of $y'$, $y$ is also comparable
to every element of $J$. Claim (a) follows.

Continuing with the proof in Case 2, in view of Claim (a), we may
suppose that $x_1$ is not meet irreducible (in $J$). Then by Lemma
\ref{lem1} (with $\mu=x_1$), there exists a unique incomparable
pair $(\alpha , \beta)$ of meet irreducibles (in $J$) such that
$x_1=\alpha\wedge\beta$. In view of the fact that $y'$ is a cover
of $x_1$, we obtain that $y'$ is equal to $A$ or $B$ ($A,B$ being
as in \S \ref{Wij}), say $y'=A$; this in turn implies that $y=C$
($C$ being as in \S \ref{Wij}; note that by hypothesis, $y$ is the
unique element covering $y'$ in $J$). Therefore we obtain that
$e_{x_1}-e_{y'}(=e_{x_1}-e_{A}),e_{y'}-e_y(=e_A-e_C)$ are in
$W\left(\tau_{\alpha,\beta}\right)$.

This completes the proof of the assertion in Case 2.

\noindent\textbf{Case 3:} Let $y'$ be both meet irreducible and
join irreducible in $J$.

If $y'$ is comparable to every other element of $J$, then $y$ is
also comparable to every other element of $J$, since by
hypothesis, $y$ is the unique element covering $ y'$ in $J$; and
the result follows.

Let then there exist a $z\in J$ such that $z$ and $y'$ are
incomparable. This in particular implies that $y'\neq
{\hat{0}}_J$; let $x\in J$ be covered by $y'$ (in fact, by
hypothesis, $x$ is unique). Proceeding as in Case 2 (especially,
the proof of Claim (a)), we obtain that $x$ is non-meet
irreducible.  Hence taking $\mu=x$ in Lemma \ref{lem1} and
proceeding as in Case 2, we obtain that $e_{y'}-e_y$ is in
$W\left(\tau_{\alpha,\beta}\right)$ ($(\alpha , \beta)$ being the
incomparable pair of irreducibles determined by $\mu$).

This completes the proof of the Lemma.
\end{proof}

\begin{theorem} Let $\tau$ be a face of $\sigma$ such that
$D_\tau$ is not contained in any $\L_{\alpha,\beta}$, for all
incomparable pair $(\alpha, \beta )$ of irreducibles in ${\L}$; in
other words $\tau$ does not contain any $\tau_{\alpha,\beta}$.
Then $\tau$ is nonsingular.
\end{theorem}
\begin{proof} As in Definition \ref{gen}, let\[W\left(\tau\right)=\{v\in W\mid f_{I_{\alpha }}
\left(v\right)=0,\,\forall\,\alpha\in D_\tau\}.\] Then
$W\left(\tau\right)$ gives a set of generators for $\tau$. By
Remark \ref{basis-sing} and \S 2.1 of \cite{F}, for $\tau$ to be
nonsingular, it must be generated by part of a basis for $N$ ($N$
being as in \S\,\ref{first-sect}).  If $W(\tau)$ is linearly
independent, then it would follow that $\tau$ is non-singular.
(Generally this is not enough to prove that $\tau$ is nonsingular;
but since all generators in $W$ have coefficients equal to $\pm
1$, any linearly independent subset of $W$ will serve as part of a
basis for $N$.)

If possible, let us assume that $W(\tau)$ is linearly dependent.
Recall that the elements of $W$ can be represented as all the line
segments in the lattice $J$, with the exception of $e_{\hat{1}}$.
Therefore, the linearly dependent generators $W(\tau)$ of $\tau$
must represent a ``loop''  of line segments in $J$. This loop will
have at least one bottom corner, left corner, top corner, and
right corner.

Let us fix an incomparable pair $(\alpha,\beta)$ of irreducibles
in ${\L}$. By Theorem \ref{main}, we have that
$W\left(\tau_{\alpha,\beta}\right)=\{e_\mu-e_A,\,e_\mu-e_B,
\,e_A-e_C,\,e_B-e_C\}$ (notation being as in that Theorem). These
four generators are represented by the four sides of a diamond in
$J$.   Thus, by hypothesis, the generators of $\tau$ represent a
loop in $J$ that does not traverse all four sides of the diamond
representing all four generators of $\tau_{\alpha,\beta}$. We have
the following identification for $\L_{\alpha,\beta}$:
$$\L_{\alpha,\beta}=\{x\in {\L}\,|\,f_{I_{x}}\equiv 0\mathrm{\,on\,}
W(\tau_{\alpha,\beta})\}.\eqno{(\dagger)}$$ The above
identification for $\L_{\alpha,\beta}$ together with the
hypothesis that $D_{\tau}\not\subseteq\L_{\alpha,\beta}$ implies
the existence of a $\theta\in D_{\tau}\cap
[\alpha\wedge\beta,\alpha\vee\beta]$; note that by $(\dagger)$, we
have
$$f_{I_{\theta}}\not\equiv 0 \mathrm{\,on\,}
W(\tau_{\alpha,\beta})$$ This implies in particular that
$\theta\not\geq C$ ($C$ being as the proof of Theorem~\ref{main});
also, $\theta\geq\mu(=\alpha\wedge\beta)$, since $\theta\in
[\alpha\wedge\beta,\alpha\vee\beta]$. Based on how $\theta$
compares to both $A$ and $B$, we can eliminate certain elements of
$W$ from $W\left(\tau\right)$. There are four possibilities; we
list all four, as well as the corresponding generators in
$W(\tau_{\alpha,\beta})$ which are not in $W(\tau)$, i.e., those
generators $v$ in $W(\tau_{\alpha,\beta})$ such that
$f_{I_{\theta }}(v)\not= 0$:
\begin{eqnarray*}
\theta\not\geq A,\,\theta\not\geq B &\Rightarrow & e_{\mu}-e_A,\,e_{\mu}-e_B\not\in W(\tau)\\
\theta\geq A,\,\theta\not\geq B & \Rightarrow &
e_A-e_C,\,e_{\mu}-e_B\not\in W(\tau)\\
\theta\not\geq A,\,\theta\geq B & \Rightarrow & e_{\mu}-e_A,\,e_B-e_C\not\in W(\tau)\\
\theta\geq A,\,\theta\geq B & \Rightarrow &
e_A-e_C,\,e_B-e_C\not\in W(\tau)
\end{eqnarray*}
Therefore, we obtain $$\mathrm{neither\,} \{e_{\mu}-e_A,e_A-e_C\}
\mathrm{\,nor\,} \{e_{\mu}-e_B,e_B-e_C\} \mathrm{\,is\,
contained\, in\,} W(\tau)\eqno{(\ast )}$$
 for any $\tau_{\alpha,\beta}$ ($(\alpha,\beta )$ being
 an incomparable pair of irreducibles
 in ${\L}$).

 Let $y',z'$ denote respectively, the left and right corners
 of our loop; let $(y,y'), (z,z')$ denote the corresponding covers
 (in $J$) which are contained in our loop. Now $y',z'$ are
 non-comparable; hence, by

 \noindent Lemma~\ref{s-3} we obtain that $(y,y')$
 (resp. $(z,z')$) are contained in some
 $W\left(\tau_{\alpha,\beta}\right)$
 (resp. $W\left(\tau_{\alpha',\beta'}\right)$). Hence we obtain
 (by Theorem \ref{main}, with notation as in that Theorem)
 $$\{e_{\mu}-e_{y'}, e_{y'}-e_{y}\}=\{e_{\mu}-e_{A}, e_{A}-e_{C}\}
 {\mathrm{\,or\,}}\{e_{\mu}-e_{B}, e_{B}-e_{C}\}$$
 But this contradicts $(\ast )$. Thus our loop
in $J$ that represented $W(\tau)$ cannot have both left and right
corners; therefore $W(\tau)$ is not a loop at all, a
contradiction. Hence, our assumption (that $W(\tau)$ is linearly
dependent) is wrong, and the result follows.
\end{proof}

Combining the above Theorem with Theorem \ref{sing} and Lemma
\ref{dim-tau-ab}, we obtain our first main Theorem:
\begin{theorem}\label{mainth} Let ${\L}$ be a distributive lattice
 such that
$J\left(\L\right)$ is a grid lattice. Then
\begin{enumerate}
 \item $\displaystyle Sing\,
X\left(\L\right)=\bigcup_{(\alpha,\beta )}
\overline{O}_{\tau_{\alpha,\beta}}$, the union being taken over
all incomparable pairs $(\alpha,\beta)$ of irreducibles in ${\L}$.
\item Sing\,$X\left(\L\right)$ is pure of codimension $3$ in
$X\left(\L\right)$; further, the generic singularities are of cone
type (more precisely, the singularity type is the same as that at the
vertex of the cone over the quadric surface $x_1x_4-x_2x_3=0$ in
$\mathbb{P}^3$ ).
\end{enumerate}
\end{theorem}

\section{Singular Loci of Bruhat-Hibi Toric Varieties}\label{mini} In this
section, we prove results for Bruhat-Hibi Toric Varieties. We
first start with recalling minuscule $G/P$'s \subsection{Minuscule
Weights and Lattices}

Let $G$ be a semisimple, simply connected algebraic group.  Let
$T$ be a maximal torus in $G$.  Let $X\left(T\right)$
 be the character
group of $T$, and $B$ a Borel subgroup containing $T$.  Let $R$ be
the root system of $G$ relative to $T$; let $R^+$ (resp.
$S=\{\alpha_1,\cdots,\alpha_l\}$) be the set of positive (resp.
simple) roots in $R$ relative to $B$ (here, $l$ is the rank of
$G$). Let $\{\omega_i, 1\le i\le l\}$ be the fundamental weights.
Let $W$ be the Weyl group of $G$, and $(,)$ a $W$-invariant inner
product on $X(T)\otimes \mathbb{R}$. For generalities on
semisimple algebraic groups, we refer the reader to \cite{B}.

Let $P$ be a maximal parabolic subgroup of $G$ with $\omega$ as
the associated fundamental weight. Let $W_P$ be the Weyl group of
$P$ (note that $W_P$ is the subgroup of $W$ generated by
$\{s_\alpha\mid \alpha\in S_P\}$). Let $W^P=W / W_P$. We have that
the Schubert varieties of $G/P$ are indexed by $W^P$, and thus
$W^P$ can be given the partial order induced by the inclusion of
Schubert varieties.

\begin{definition} A fundamental weight $\omega$ is called
\emph{minuscule} if
$\left<\omega,\beta\right>(={\frac{2(\omega,\beta)}{(\beta,\beta)}})\leq
1$ for all $\beta\in R^+$; the maximal parabolic subgroup
associated to $\omega$ is called a \emph{minuscule parabolic
subgroup}.
\end{definition}

\begin{remark}[cf \cite{Hiller}] Let $P$ be a maximal parabolic
subgroup; if $P$ is minuscule then $W/W_P$ is a distributive
lattice.
\end{remark}

\begin{definition} For $P$ a minuscule parabolic subgroup,
we call $\mathcal{L}=W / W_P$ a \emph{minuscule lattice}.
\end{definition}
\begin{definition} We call $X({\L})$ a
\emph{Bruhat-Hibi toric variety} (B-H toric variety for short) if
${\L}$ is a minuscule lattice.
\end{definition}

In order to begin work on these B-H toric varieties, we first list
all of the minuscule fundamental weights.  Following the indexing
of the simple roots as in \cite{bou}, we have the complete list of
minuscule weights for each type:
\[ \begin{array}{lll}
        \mbox{Type } \mathbf{A_n} & : & \mbox{Every fundamental weight is minuscule}\\
        \mbox{Type } \mathbf{B_n} & : & \omega_n\\
        \mbox{Type } \mathbf{C_n} & : & \omega_1\\
        \mbox{Type } \mathbf{D_n} & : & \omega_1,\omega_{n-1},\omega_n\\
        \mbox{Type } \mathbf{E_6} & : & \omega_1, \omega_6\\
        \mbox{Type } \mathbf{E_7} & : & \omega_7.
\end{array}\]
There are no minuscule weights in types $\mathbf{E_8},
\mathbf{F_4}, \mbox{ or } \mathbf{G_2}.$

Before proving that each minuscule lattice has grid lattice join
irreducibles, we must introduce some additional lattice notation.
For a poset $P$, let $\mathcal{I}\left(P\right)$ represent the
lattice of ideals of $P$.  Thus for a distributive lattice $\L$,
$\L=\mathcal{I}\left( J\left(\L\right)\right)$ (cf. Theorem
\ref{5.10}).  (Notice that the empty set is considered the minimal
ideal, and in Theorem \ref{5.10} we do not include the minimal
element in $P$. Therefore, in this section,
$\mathcal{I}\left(J\right)$ will have a minimal element that is
not an element of $J$.)

For $k\in\N$, let $\underline{k} $ be the the totally ordered set
with $k$ elements.  The symbols $\oplus$ and $\times$ denote the
disjoint union and (Cartesian) product of posets.

Let $\mathbf{X_n}\left(\omega_i\right)$ denote the minuscule
lattice $W/W_P$ where $P$ is a parabolic subgroup associated to
$\omega_i$ in the root system of type $\mathbf{X_n}$.

\begin{theorem} [cf. \cite{Proctor} Propositions 3.2 and 4.1\footnote{Our notation differs significantly than that used in \cite{Proctor}; namely that where we use $\mathcal{I}$, Proctor uses $J$; whereas we use $J$ to signify the set of join irreducibles.}]\label{proctor} The minuscule lattices have the following combinatorial descriptions:
\[\begin{array}{l}
\mathbf{A_{n-1}}\left(\omega_j\right)\cong \mathcal{I}\left(\mathcal{I}\left(\underline{j-1}\oplus \underline{n-j-1}\right)\right) \\
\mathbf{C_{n}}\left(\omega_1\right)\cong \underline{2n}\\
\mathbf{B_{n}}\left(\omega_n\right)\cong \mathbf{D_{n+1}}\left(\omega_{n+1}\right)\cong\mathbf{D_{n+1}}\left(\omega_n\right)\cong\mathcal{I}\left( \mathcal{I}\left(\mathcal{I}\left(\underline{1}\oplus \underline{n-2}\right)\right)\right) \\
\mathbf{D_{n}}\left(\omega_1\right)\cong \mathcal{I} ^{n-1}\left(\underline{1}\oplus \underline{1} \right) \\
\mathbf{E_{6}}\left(\omega_1\right)\cong\mathbf{E_6}\left(\omega_6\right)\cong \mathcal{I} ^4 \left(\underline{1}\oplus \underline{2}\right) \\
\mathbf{E_{7}}\left(\omega_7\right)\cong \mathcal{I}
^5\left(\underline{1}\oplus \underline{2}\right) .
\end{array}\]
\end{theorem}

This theorem is very convenient in working with the faces of B-H
toric varieties, because the join irreducible lattice of each of
these minuscule lattices is very easy to see, simply by
eliminating one $\mathcal{I}\left(\cdot\right)$ operation.  Our
goal is to show that each minuscule lattice has join irreducibles
with grid lattice structure.

\subsection{Minuscule lattices $\mathbf{A_{n-1}}\left(\omega_j\right)$}

\begin{remark}[cf. \cite{Proctor} Proposition 4.2] The join irreducibles of the minuscule lattice $\mathbf{A_{n-1}}\left(\omega_j\right)$ are isomorphic to the lattice
\[ \underline{j}\times\underline{n-j}.\]
\end{remark}

Therefore, every element of
$J\left(\mathbf{A_{n-1}}\left(\omega_j\right)\right)$ can be
written as the pair $(a,b)$, for $1\leq a\leq j ,\, 1\leq b\leq
n-j$.  This leads us to the following result,

\begin{cor} The minuscule lattice $\mathbf{A_{n-1}}\left(\omega_j\right)$ has grid lattice join irreducibles.
\end{cor}

Note that the result about the singular loci of B-H toric
varieties of type $\mathbf{A_n}\left(\omega_j\right)$ was already
proved in \cite{G-H-paper}, as well as more results about the
multiplicities of singular points, but using the unique
combinatorics of these lattices.

\subsection{Minuscule lattices $\mathbf{C_{n}}\left(\omega_1\right)$
} This minuscule lattice is totally ordered, and the associated
B-H toric variety is simply the affine space of dimension $2n$.

\subsection{Minuscule lattices $\mathbf{B_{n-1}}
\left(\omega_{n-1}\right)\cong
\mathbf{D_{n}}\left(\omega_{n-1}\right)\cong\mathbf{D_{n}}\left(\omega_n\right)$}

From Theorem \ref{proctor}, we have
\[J\left(\mathbf{D_n}\left(\omega_n\right)\right)\cong \mathcal{I} ^2\left(\underline{1}\oplus \underline{n-3}\right)\cong \mathbf{A_{n-1}}\left(\omega_2\right).\]
It is a well known result that $\mathbf{A_{n-1}}\left(2\right)$
represents the lattice of Schubert varieties in the Grassmannian
of 2-planes in $K^{n}$, and the Schubert varieties are indexed by
$I_{2,n}=\{(i_1,i_2 )\mid 1\leq i_1 <i_2\leq n\}$. Therefore,
\[J\left(\mathbf{D_n}\left(\omega_n\right)\right)\cong I_{2,n}.\]
The lattice $I_{2,n}$ is therefore distributive, (being another
minuscule lattice), and clearly a grid lattice.  This leads to the
following result,

\begin{cor} The minuscule lattices
$\mathbf{B_{n-1}}\left(\omega_{n-1}\right),\,
\mathbf{D_{n}}\left(\omega_{n-1}\right)$,

\noindent $\mathbf{D_{n}}\left(\omega _n\right)$ have grid lattice
join irreducibles.
\end{cor}

\subsection{Minuscule lattices $\mathbf{D_n}\left(\omega_1\right)$}

From Theorem \ref{proctor}, we have

\noindent $J\left(\mathbf{D_n}\left(\omega_1\right)\right)\cong
\mathcal{I} ^{n-2}\left(\underline{1}\oplus\underline{1}\right)$.
This lattice of join irreducibles is isomorphic to the following sublattice of
$\N\times\N$ (drawn horizontally):
\[ \xymatrix@-10pt{
& & (2,n-2)\ar@{-}[r] & (2,n-1)\ar@{--}[rr] & & (n,n-1)\\
 (1,1)\ar@{--}[rr] & & (1,n-2)\ar@{-}[u]\ar@{-}[r] & (1,n-1)\ar@{-}[u] & & }\] 
Clearly this is a grid lattice.

\subsection{Minuscule lattices $\mathbf{E_6}\left(\omega_1\right)\cong\mathbf{E_6}\left(\omega_6\right)$, and $\mathbf{E_7}\left(\omega_7\right)$}

Let
$\mathbf{H_6}=\mathbf{E_6}\left(\omega_1\right)=\mathbf{E_6}\left(\omega_6\right)$
and $\mathbf{H_7}=\mathbf{E_7}\left(\omega_7\right)$.  Since there
are only two exceptional cases, it is best to explicitly give the
grid lattice structure to the join irreducibles. Thus, we have the
two join irreducible lattices below, with each lattice point given
coordinates in $\N\times\N$.  Coincidentally,
$J\left(\mathbf{H_6}\right)=\mathbf{D_5}\left(\omega_5\right)$ and
$J\left(\mathbf{H_7}\right)=\mathbf{H_6}$.

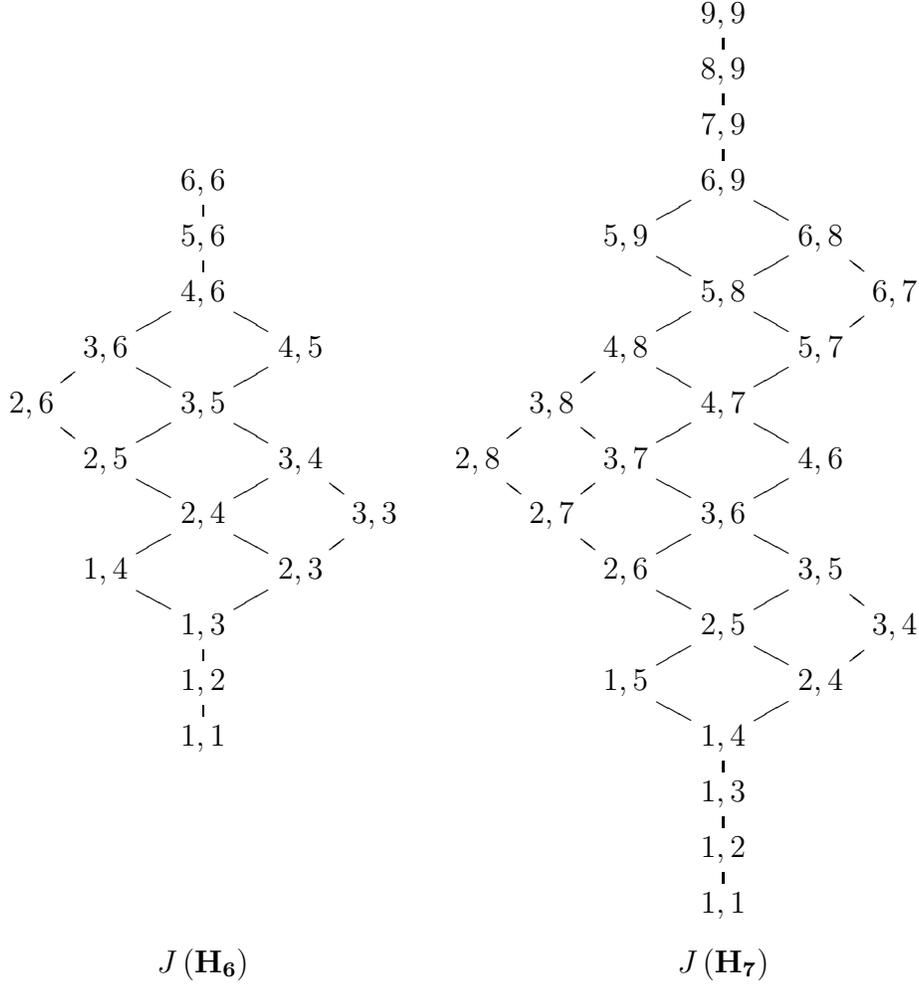
\[ \xymatrix@-20pt{
& & & & & & & & & 9,9 \ar@{-}[d]& & \\
& & & & & & & & & 8,9\ar@{-}[d] & & \\
& & & & & & & & & 7,9\ar@{-}[d] & & \\
& & 6,6 \ar@{-}[d]& & & & & & & 6,9\ar@{-}[dl]\ar@{-}[dr] & & \\
& & 5,6 \ar@{-}[d]& & & & & & 5,9\ar@{-}[dr] & & 6,8\ar@{-}[dr]\ar@{-}[dl] & \\
& & 4,6\ar@{-}[dr]\ar@{-}[dl]& & & & & & & 5,8\ar@{-}[dr]\ar@{-}[dl] & & 6,7\ar@{-}[dl] \\
& 3,6 \ar@{-}[dr]\ar@{-}[dl]& & 4,5\ar@{-}[dl] & & & & & 4,8\ar@{-}[dr]\ar@{-}[dl]& & 5,7\ar@{-}[dl]& \\
2,6 \ar@{-}[dr]& & 3,5\ar@{-}[dr]\ar@{-}[dl] & & & & & 3,8\ar@{-}[dr]\ar@{-}[dl] & & 4,7\ar@{-}[dr]\ar@{-}[dl] & & \\
& 2,5\ar@{-}[dr]& & 3,4\ar@{-}[dr]\ar@{-}[dl] & & & 2,8\ar@{-}[dr] & & 3,7\ar@{-}[dr]\ar@{-}[dl] & & 4,6\ar@{-}[dl] & \\
& & 2,4 \ar@{-}[dr]\ar@{-}[dl]& & 3,3\ar@{-}[dl] & & & 2,7\ar@{-}[dr] & & 3,6\ar@{-}[dr]\ar@{-}[dl] & & \\
& 1,4\ar@{-}[dr]& & 2,3\ar@{-}[dl] & & & & & 2,6 \ar@{-}[dr]& & 3,5\ar@{-}[dr]\ar@{-}[dl] & \\
& & 1,3\ar@{-}[d]& & & & & & & 2,5\ar@{-}[dr]\ar@{-}[dl]& & 3,4\ar@{-}[dl]\\
& & 1,2\ar@{-}[d]& & & & & & 1,5\ar@{-}[dr] & & 2,4\ar@{-}[dl] & \\
& & 1,1 & & & & & & & 1,4 \ar@{-}[d]& & \\
& & & & & & & & & 1,3 \ar@{-}[d]& & \\
& & & & & & & & & 1,2 \ar@{-}[d]& & \\
& & & & & & & & & 1,1 & & \\
& & J\left(\mathbf{H_6}\right) & & & & & & &
J\left(\mathbf{H_7}\right)& &
 } \]

This completes the individual discussion for each type of
minuscule lattice, leading us to the following result.

\begin{cor} If $\L$ is a minuscule lattice, then $J\left(\L\right)$ is a grid lattice.
\end{cor}

\noindent Thus, for $\L$ any minuscule lattice, letting
\[\Phi = \{(\alpha,\beta )\mid \alpha,\beta\mbox{ non-comparable
 irreducibles in }\L\},\]
we have completed the proof of the conjecture from
\cite{G-H-paper}, thanks to Theorem \ref{mainth}:

\begin{theorem}\label{theth} For the B-H toric variety $X\left(\L\right)$,
\[ Sing\, X\left(\L\right) = \bigcup _{(\alpha,\beta )\in\Phi}\overline{O}_{\tau_{\alpha,\beta}}.\]
In other words, $X\left(\L\right)$ is smooth at $P_\tau$ ($\tau$
being a face of $\sigma$) if and only if for each pair
$(\alpha,\beta )\in\Phi$, there exists at least one $\gamma\in
[\alpha\wedge\beta,\alpha\vee\beta]$ such that $P_\tau(\gamma)$ is
non-zero.
\end{theorem}



\begin{thebibliography}{99}


\bibitem{Font}
Victor V. Batyrev, Ionut Ciocan-Fontanine, Bumsig Kim, Duco van
Straten, {\em Mirror Symmetry and Toric Degenerations of Partial
Flag Manifolds}, Acta Math. 184 (2000), no. 1, 1--39.

\bibitem{Font2}
Victor V. Batyrev, Ionut Ciocan-Fontanine, Bumsig Kim, Duco van
Straten, {\em Conifold Transitions and Mirror Symmetry for
Calabi-Yau Complete Intersections in Grassmannians}, Nuclear Phys.
B 514 (1998), no. 3, 640--666.

\bibitem{B}
A. Borel, {\em Linear Algebraic Groups}, second edition,
Springer-Verlag, New York, 1991.

\bibitem{bou}
{N. Bourbaki},  Groupes et Alg\` ebres de Lie, Chapitres 4, 5 et
6, {\em Hermann}, Paris 1968.

\bibitem{G-H-paper}
J. Brown and V. Lakshmibai, {\em Singular Loci of Grassmann-Hibi
Toric Varieties}, submitted to Adv. Math.

\bibitem{Eisenbud}
D. Eisenbud, {\em Commutative Algebra with a View Toward Algebraic
Geometry}, Springer-Verlag, New York, 1995.

\bibitem{ES}
D. Eisenbud and B. Sturmfels, {\em Binomial ideals}, Duke
Mathematical Journal \textbf{84} (1996), 1--45.

\bibitem{F}
W. Fulton, {\em Introduction to Toric Varieties}, Annals of Math.
Studies {\bf 131}, Princeton U. P., Princeton N. J., 1993.

\bibitem{GLdef}
N. Gonciulea and V. Lakshmibai, {\em Degenerations of flag and
Schubert varieties to toric varieties}, Transformation Groups, vol
1, no:3 (1996), 215-248.

\bibitem{g-l}
N. Gonciulea and V. Lakshmibai, {\em Schubert varieties, toric
varieties and ladder determinantal varieties }, Ann. Inst.
Fourier, t.47, (1997), 1013-1064.

\bibitem{Ha}
R. Hartshorne, {\em Algebraic Geometry}, Springer-Verlag, New
York, 1977.

\bibitem{Hi}
T. Hibi, {\em Distributive lattices, affine semigroup rings, and
algebras with straightening laws}, Commutative Algebra and
Combinatorics, Advanced Studies in Pure Math. {\bf 11} (1987)
93-109.

\bibitem{Hiller}
H. Hiller, {\em Geometry of Coxeter groups}, Pitman Adv. pub.
prog., Research Notes in Math, vol 54, Boston, London, Melbourne,
1982

\bibitem{KK} G. Kempf et al, {\em Toroidal Embeddings}, Lecture notes in
Mathematics, N0. 339,  Springer-Verlag, 1973.

\bibitem{L-G} V. Lakshmibai and N. Gonciulea, {\em Flag Varieties}, Hermann, \`Editeurs des Sciences et des Arts, 2001.

\bibitem{LM}
V. Lakshmibai and H. Mukherjee, {\em Singular Loci of Hibi Toric
Varieties}, submitted to J. Ramanujan Math. Soc.

\bibitem{l-mu2}
V. Lakshmibai and H. Mukherjee, {\em Standard Monomial Basis for
Tangent Cones at Singular Points of Hibi Toric Varieties},
preprint.

\bibitem{Proctor}
R.A. Proctor, {\em Bruhat Lattices, Plane Partition Generating
Functions, and Minuscule Representations}, European J.
Combinatorics {\bf 5} (1984), 331-350.


\bibitem{St}
B. Sturmfels, {\em Gr\" obner bases and convex polytopes},
American Mathematical Society, Univ. Lecture Series, No. 8,
Providence, Rhode Island, 2002.

\bibitem{Wa}
D. G. Wagner, {\em Singularities of toric varieties associated
with finite distributive lattices}, Journal of Algebraic
Combinatorics \textbf{5} (1996), 149--165.

\end{thebibliography}
\end{document}